\documentclass[a4paper,12pt]{amsart}

\usepackage{verbatim}
\usepackage{amsfonts,amsmath,amssymb}
\usepackage{amssymb}
\usepackage[dvips]{graphicx}
\usepackage{psfrag}
\usepackage{epsfig}

\usepackage{color}
\usepackage{graphicx}
\usepackage{graphics}

\mathchardef\emptyset="001F

\theoremstyle{plain}

\newtheorem{theorem}{Theorem}[section]
\newtheorem{lemma}[theorem]{Lemma}
\newtheorem{proposition}[theorem]{Proposition}

\newtheorem{remark}[theorem]{Remark}

\theoremstyle{definition}
\theoremstyle{remark}
\numberwithin{equation}{section}

\newcommand{\e}{\varepsilon}
\newcommand{\Om}{\Omega}

\newcommand{\weak}{\rightharpoonup}

\newcommand{\R}{{\mathbb R}}
\newcommand{\F}{{\mathcal F}}
\newcommand{\G}{{\mathcal G}}
\newcommand{\C}{{\mathcal C}}
\newcommand{\I}{{\mathcal I}}
\newcommand{\E}{{\mathcal E}}
\renewcommand{\O}{{\mathcal O}}

\newcommand{\Z}{{\mathbb Z}}

\newcommand{\M}{{\mathbb M}}
\renewcommand{\H}{{\mathcal H}}

\newcommand{\Mtt}{\M^{3\times 3}}

\newcommand{\var}{\varphi}

\newcommand{\no}{\noindent}
\newcommand{\non}{\nonumber}

\newcommand{\diag}{\hbox{{\rm diag}}}
\newcommand{\dsp}{\displaystyle}

\newcommand{\dist}{\mathrm{dist}}
\newcommand{\ul}{u_{\lambda}}
\newcommand{\Fh}{F^{(h)}}
\newcommand{\Fhh}{F^{(h)}_{h}}
\newcommand{\gh}{\gamma_{_{H}}}
\newcommand{\ghh}{\hat\gamma_{_{H}}}
\newcommand{\vh}{v^{(h)}}
\newcommand{\uh}{u^{(h)}}
\newcommand{\tv}{\tilde{v}}
\newcommand{\curl}{{\rm curl}\,}
\renewcommand{\H}{{\mathcal H}}
\renewcommand{\b}{{\bf b}}
\newcommand{\co}{{\rm co}}

\newcommand{\tp}{\tilde{p}}

\title[]{Derivation of a rod theory for biphase materials with dislocations 
at the interface}

\author[Stefan M\"uller]{Stefan M\"uller$^{1}$}
\author[Mariapia Palombaro]{Mariapia Palombaro$^{2}$}

\begin{document}
\baselineskip3.15ex
\vskip .3truecm

\maketitle
\small
\noindent
$^1$ Hausdorff Center for Mathematics $\&$ Institute for Applied Mathematics, Universit\"at Bonn, Endenicher Allee 60, 53115 Bonn, Germany.
Email: sm@hcm.uni-bonn.de\\
\noindent
$^2$ SISSA,Via Beirut 2-4, 34014 Trieste, Italy.
Email:palombar@sissa.it

\vspace{2mm}

\begin{abstract}
Starting from three-dimensional elasticity we derive a 
rod theory for biphase materials with a prescribed dislocation at 
the interface. 
The stored energy density is assumed to be non-negative and to vanish 
on a set consisting of two copies of $SO(3)$. 
First, we rigorously justify the assumption of dislocations at the interface.
Then, we consider the typical scaling 
of multiphase materials and we perform an asymptotic study of the rescaled energy, 
as the diameter of the rod goes to zero, in the framework of $\Gamma$-convergence.
 

\end{abstract}

\section{Introduction}
\no
We study the behavior of an elastic thin beam consisting of two parts made of 
different materials. The interface between the two parts of the beam is fixed. 
Our objective is, first, to rigorously prove that 
formation of dislocations on such interface is energetically 
more favorable than purely elastic deformation when the radius of the 
cross-section is sufficiently large. 
Second, to derive a one-dimensional theory of elastic thin beams with a 
prescribed dislocation on the interface.

The motivation to look at this problem relies on the connection with 
the study of nanowire heterostructures, which have  
important applications in semiconductor electronics. 
A heterostructure is a material obtained through an epitaxial growth process, 
where two materials featuring different lattice constants are brought together 
by deposition of one material (the overlayer) on top of the other (the underlayer). 
In general, lattice mismatch will prevent growth of defect-free epitaxial film over a 
substrate unless the thickness of the film is below certain critical thickness; in this 
last case lattice mismatch is compensated by the strain in the film. 
In contrast, as confirmed by experimental observations, one-dimensional 
systems, i.e., longitudinally heterostructured nanowires, can be grown defect-free 
more readily than their two-dimensional counterparts. A better understanding of 
nanowires is therefore crucial in the study and use of heterostructures.


A schematic of a heterostructured nanowire is showed in Figure 1. The radii of 
the unstrained underlayer and overlayer are denoted by $R$ and $r$ respectively. 
The lattice mismatch, $\alpha$, between the overlayer and the underlayer is defined 
as

\begin{equation}\label{mismatch}
\alpha:=1-\frac{r}{R}\,.
\end{equation}

\no
For a given mismatch $\alpha$, if the radii $R$ and $r$ are sufficiently small, the system 
is elastically strained and no dislocation arises. Ultimately, as the radii increase, 
the mismatch strain is relieved by formation of misfit dislocations at the interface. 
In the dislocated system, a small portion of the total mismatch $\alpha$ is accomodated 
by the dislocations, while the remainder (the residual mismatch) is accomodated by elastic 
strain both in the underlayer and overlayer (see Figure 1). 
Figure 2 represents a longitudinal section of a dislocated nanowire in the atomistic picture (where the crystalline lattice is assumed to be cubic): 
we observe an additional row of atoms in the overlayer.

A model for the critical radius for which the first dislocation appears has been developed, 
e.g., in \cite{egcs} in the context of linearized elasticity. The critical radius 
$R$ is described in \cite{egcs} as a function of the mismatch $\alpha$, and is shown to be 
roughly an order of magnitude larger than the critical thickness of the corresponding thin 
film/substrate system.


\begin{figure}
\centering
\psfrag{R}{$R$}
\psfrag{r}{$r$}
\psfrag{d}{\tiny{dislocation line}}
\includegraphics[scale=0.9]{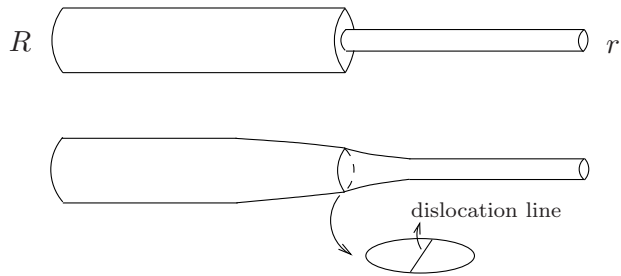}
\caption{Schematic of a nanowire heterostructure before and after 
interfacial bonding.}
\end{figure}


\begin{figure}
\centering
\includegraphics[scale=0.9]{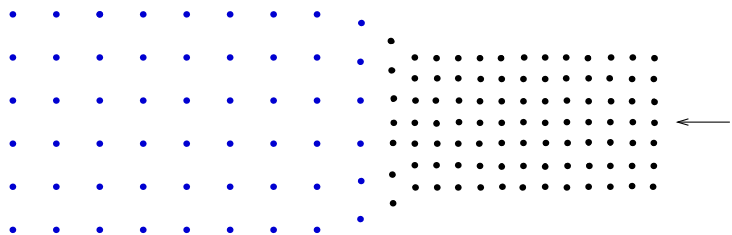}
\caption{Longitudinal section of the beam in the atomistic picture.}
\end{figure}


The purpose of this paper is, first, to rigorously justify formation of dislocations on the 
interface between the two parts of the beam; then, to derive a one dimensional model 
describing the deformations of the beam with a given dislocation.
For the second part we consider the case of one misfit dislocation, 
though our analysis extends  
as well to the case of more dislocations without any additional difficulty.
 
More precisely, we consider a cylindrical region $\Om_{h}:=(-L,L)\times hS$, which 
represents the reference configuration of the beam, where $S$ is the disk of radius $r$ 
in  $\R^{2}$, i.e., $S=\{x_{1}^{2}+x_{2}^{2}<r^{2}\}$, 
and $h$ is a small 
positive parameter, which, in the atomistic picture, is of the order of the atomic distance. 
Theorem \ref{mainthm35} shows that when $r$ is sufficiently large, formation of dislocations 
is energetically more favorable than purely elastic deformation. 

In the second part of the paper, we prescribe the dislocation.
We assume that the dislocation line, $\Gamma_{h}$, has the form 
$$
\Gamma_{h}:=h\Gamma\,,\quad\quad \Gamma\subset S\,,
$$
\no
where $\Gamma$ is a Lipschitz, relatively closed curve in $S$. 
The latter condition implies that 
$\Om_h\setminus\Gamma_h$ is not simply connected. 
 
We assume that the elastic energy (per unit cross-section) has the form

\begin{equation}\label{startenergy}
\E^{(h)}(G):=\frac{1}{h^{2}}\int_{\Om_{h}}W(z,G(z))\,dz \,,
\end{equation}

\no
where $\dsp G\in  L^{p}(\Om_{h};\Mtt)$ satisfies:

\begin{equation}\label{disloconstraint}
\curl G = -h\b\otimes \dot\Gamma_{h} \,d\H^{1}\llcorner \Gamma_{h}
\end{equation}

\begin{figure}
\centering
\psfrag{x}{$x_1$}
\psfrag{y}{$x_2$}
\psfrag{z}{$x_3$}
\psfrag{h}{\small{$2\hspace{-0.2mm}hr$}}
\psfrag{Gammah}{$\Gamma_h$}
\psfrag{-L}{$-L$}
\psfrag{L}{$L$}
\includegraphics[scale=0.9]{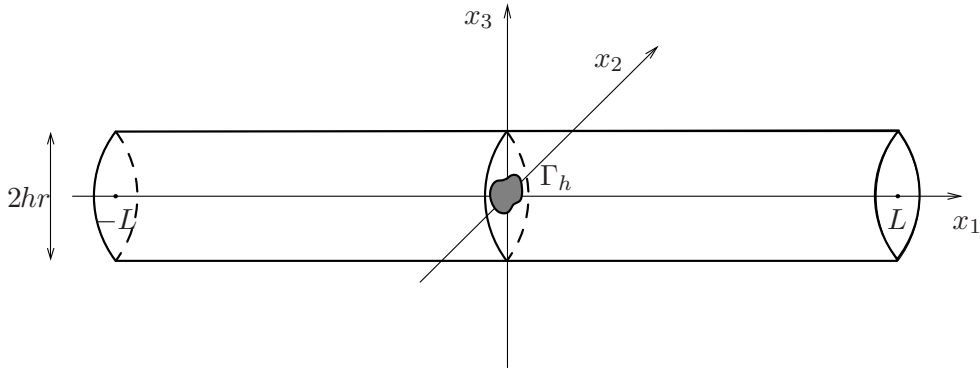}
\caption{Reference configuration of the beam.}
\label{reference}
\end{figure}

\no
in the sense of distributions. In \eqref{disloconstraint}, the vector $h\b$, with 
$\b\in \R^{3}, |\b|=1$,    
denotes the Burgers vector, which, together with the dislocation line, uniquely 
characterizes the dislocation. 
We observe that any field $G$ satisfying \eqref{disloconstraint} is locally the gradient of a 
Sobolev map. More precisely, if $\omega\subset\Om_{h}\setminus\Gamma_{h}$ is simply 
connected, then there exists $u\in W^{1,p}(\omega;\R^{3})$ such that 
$G=Du$ a.e. in $\omega$. 
In particular, if $\Gamma$ is a closed loop in $S$, 
one can 
take $\omega=\Om_{h}\setminus D_{h}$, where $D_{h}:=hD$, 
and $D$ is the flat region enclosed by the curve $\Gamma$ 
($D_{h}$ is the shadowed set in Figure \ref{reference}).
Then, $G=\nabla u$ a.e. in $\Om_{h}$, 
where $u\in SBV(\Om_{h};\R^{3})$ and  its distributional gradient satisfies

\begin{equation*} 
Du=\nabla u \,dx + h\b\otimes e_{1} \,d\H^{2}\llcorner D_{h} \,.
\end{equation*}
 
\no
Therefore $G=\nabla u$ is the absolutely continuous part (with respect to Lebesgue measure) 
of the gradient $Du$. Following \cite{ortiz-vienna}, we interpret $G$ as 
the elastic part of the deformation. 
More in general, $G$ may be regarded as the elastic part of a deformation which has a 
constant jump, equal to $h\b$, across any surface having $\Gamma_{h}$ as its boundary. 

The domain of the energy functional \eqref{startenergy} is thus defined as 

\begin{equation*}
\G^{(h)}:=\{G\in  L^{p}(\Om_{h};\Mtt) \::\:
\curl G = -h\b\otimes \dot\Gamma_{h} \,d\H^{1}\llcorner \Gamma_{h} \}\,,
\end{equation*}

\no
where $p<2$. Indeed, because of \eqref{disloconstraint}, the fields of $\G^{(h)}$ 
cannot expect to be in $L^{2}(\Om_{h};\Mtt)$.


\no
Furthermore, we assume that the density of energy $W:\Om_{h}\to [0,+\infty)$ 
has the form

$$
W(x,A)=
\begin{cases}
W_{1}(A) & \text{ if } x_{1}\in(-L,0) \\
W_{2}(A) & \text{ if } x_{1}\in(0,L) \\
\end{cases}
$$
\no
where the functions $W_{1}$ and $W_{2}$ satisfy the following conditions \\

\no
(i) $W_{i}\in C^{0}(\Mtt)$, $i=1,2$;

\no
(ii) $W_{i}$ is frame indifferent, i.e., $W_{i}(A)=W_{i}(RA)$ for every $A\in\Mtt$ and 
$R\in SO(3)$, $i=1,2$;

\no
(iii) there exist $C_{1},C_{2}>0$, $H\in\Mtt$, with $\det H>0$,  
such that for every $A\in\Mtt$

$$
C_{1}\Big(\dist^{2}((A,SO(3)))\wedge (|A|^{p}+1)\Big)
\leq W_{1}(A)
\leq C_{2}\Big(\dist^{2}((A,SO(3)))\wedge (|A|^{p}+1)\Big)\,,
$$ 

\no
and
$$
C_{1}\Big(\dist^{2}((A,SO(3)H))\wedge (|A|^{p}+1)\Big)
\leq W_{2}(A)
\leq C_{2}\Big(\dist^{2}((A,SO(3)H))\wedge (|A|^{p}+1)\Big)\,,
$$ 

\no
for some $p\in(1,2)$.
(Remark that a typical $H$ is, for example, $H=(1-\alpha)I$, where $I$ is the identity 
matrix and $\alpha$ is defined by \eqref{mismatch}.)
More generally in the following we assume that 
$H=\diag(\zeta_{1},\zeta_{2},\zeta_{3})$, 
with $\zeta_{i}>0$ for $i=1,2,3$. 

In order to recast the functionals over varying domains $\Om_{h}$ 
into functionals with a fixed domain $\Om$, we introduce 
in \eqref{startenergy} the following 
change of variables:

$$
z_{1}=x_{1}\,, \quad z_{2}=hx_{2}\,, \quad
z_{3}=hx_{3}\,,
$$
\no
and rescale the elements of $\G^{(h)}$ accordingly

\begin{equation}\label{rescaledfield}
G(z(x))=:F_{h}(x)\,.
\end{equation}

\no
In \eqref{rescaledfield} we used the notation 

$$
F_{h}:=\left(F^{1}\Big|\,\frac{1}{h}F^{2}\Big|\,
\frac{1}{h}F^{3}\right) \,,
$$
\no
where $F^{i}$ stands for the $i$th column of $F$. 
We now rewrite \eqref{startenergy} in terms of maps 
from $\dsp\Om:=(-L,L)\times S$ to $\Mtt$:

$$
\I^{(h)}(F):=\int_{\Om}W(x,F_{h}(x))\,dx =\E^{(h)}(G)\,.
$$ 

\no
It will be convenient to define the set of admissible deformations in the fixed domain 
$\Om$: 

$$
\F^{(h)}:=\{F\in  L^{p}(\Om;\Mtt) \::\:
\curl F = -h\b\otimes \dot\Gamma \,d\H^{1}\llcorner \Gamma \}\,.
$$


\no
Our goal is to study the asymptotic behavior of the rescaled 
functionals $\frac{1}{h}\I^{(h)}(F)$ in the framework of $\Gamma$-convergence. 
This problem was already addressed in \cite{mm} 
in the dislocation-free case. The main difference here is that, due to the presence of a 
dislocation, one has to work with growth conditions slower than quadratic, as 
specified in (iii), which require suitable modifications of the methods introduced in \cite{mm}. 

In Theorem \ref{thm1} we show that if a sequence $\dsp \{\Fh\}\subset\F^{(h)}$ is such 
that $\frac{1}{h}\I^{(h)}(\Fh)<C$, then, up to subsequences, the sequence 
$\{\Fh_{h}\}$ converges weakly in 
$L^{p}(\Om;\Mtt)$ to some limit $F\in L^{p}((-L,L);\Mtt)$ such that 

\begin{equation*}
F\in 
\begin{cases}
{\rm co}(SO(3))     &  \text{a.e in }(-L,0)\,,\\
{\rm co}(SO(3)H)   &  \text{a.e in }(0,L)\,, 
\end{cases}
\end{equation*}

\no
where co$(A)$ denotes the convex hull of $A$ for any $A\in\Mtt$. 
Finally, in Theorem \ref{thm3}, we compute the $\Gamma$-limit of the sequence 
$\frac{1}{h}\I^{(h)}$.


\section{Preliminary results}

\no
Throughout this paper the letter $C$ denotes various positive constants whose 
precise value may change from place to place. Its dependence on other variables will be 
emphasized only if necessary.

We will use the following two results from \cite{fjm}.

\begin{theorem}
Let $n\geq 2$, and let $1\leq p < \infty$.
Suppose that $U\subset\R^{n}$ is a bounded Lipschitz domain. 
There exists a constant $C(U)$
such that for each $u\in W^{1,p}(U,\R^{n})$ 
there exists $R\in SO(n)$ such that 

\begin{equation}\label{rigidity}
\|D u-R\|_{L^{p}(U)} \leq C
\|\dist(D u,SO(n))\|_{L^{p}(U)}\,.
\end{equation}

\end{theorem}

\vspace{2mm}

\begin{proposition}\label{propolip} 
Let $n,m\geq 1$, and let $1\leq p < \infty$. Suppose that $U\subset\R^{n}$ 
is a bounded Lipschitz domain. Then there exists a constant $C(U, m, p)$ 
such that for each $u\in W^{1,p}(U,\R^{m})$ and each $\lambda>0$, 
there exists $u_{\lambda}:U\to\R^{m}$ such that 
\vspace{3mm}

\begin{eqnarray}
&\label{con1}\dsp \|D u_{\lambda}\|_{L^{\infty}(U)}\leq C\lambda\,,\hspace{2.7cm} \\
&\label{con2}\dsp |\{u\neq u_{\lambda}\}|\leq \frac{C}{\lambda^{p}}
\int_{\{|D u>\lambda|\}} |D u|^{p} \,dx\,, \\
&\label{con3}\hspace{1cm}\dsp\|D u-D u_{\lambda}\|^{p}_{L^{p}(U)}\leq C
\int_{\{|D u|>\lambda\}} |D u|^{p} \,dx \,.
\end{eqnarray}
\end{proposition}

The next proposition provides a generalization of the rigidity estimate \eqref{rigidity}, 
which cannot be applied as it is, due to the growth condition (iii) required for the 
function $W$.  It will be used to prove the compactness of sequences with equibounded 
energy.

\begin{proposition}\label{rigext}
Let $n\geq 2$, and let $1\leq p < 2$. Suppose that $U\subset\R^{n}$ 
is a bounded Lipschitz domain. Then there exists a constant $C(U)$ 
such that for each $u\in W^{1,p}(U,\R^{n})$ 
there exists $R\in SO(n)$ such that 


\begin{equation}\label{prigidity}
\int_{U} |D u-R|^{2}\wedge (|D u|^{p}+1) \,dx \leq C(U)
\int_{U} \dist^{2}(D u,SO(n))\wedge (|D u|^{p}+1) \,dx \,.
\end{equation}
\end{proposition}

\begin{proof}
Let $\lambda>0$ and let $\ul$ be given by Proposition  \ref{propolip}.  
Set $U_{\lambda}:=\{u=\ul\}$.
The rigidity estimate \eqref{rigidity} implies that there exists $R\in SO(n)$ 
such that 

\begin{align}\label{step1}
\non   \int_{U} |D \ul-R|^{2}\,dx & \leq 
          C \int_{U} \dist^{2}(D \ul,SO(n))\,dx \\
       & = C \int_{U_{\lambda}} \dist^{2}(Du,SO(n))\,dx +
          C \int_{U\setminus U_{\lambda}} \dist^{2}(D \ul,SO(n))\,dx \,.
\end{align}
Since $|D \ul|\leq C\lambda$, we can find a constant $C$,  
depending on $\lambda$, such that 

\begin{equation}\label{step2}
\int_{U_{\lambda}} \dist^{2}(Du,SO(n))\,dx \leq
C \int_{U_{\lambda}} \dist^{2}(D u,SO(n)) 
\wedge (|D u|^{p}+1) \,dx\,.
\end{equation}

\no
For the second term of \eqref{step1} we use \eqref{con1}-\eqref{con2}
to get, for sufficiently large $\lambda$,  

\begin{align}\label{step3}
\non  \int_{U\setminus U_{\lambda}} \dist^{2}(D \ul,SO(n))\,dx & \leq  
            C\lambda^{2}|U\setminus U_{\lambda}|\\
\non  &\leq  C\lambda^{2-p}\int_{\{|D u|>\lambda\}} |D u|^{p}\,dx \\
        & \leq C\lambda^{2-p}
           \int_{U} \dist^{2}(D u,SO(n)) \wedge (|D u|^{p}+1) \,dx \,.
\end{align}

\no
In the last inequality of \eqref{step3} we used the fact that, 
for sufficiently large $\lambda$,   
$|D u|^{p}+1<\dist^{2}(D u,SO(n))$ a.e. 
in the set $\{|D u|>\lambda\}$.  
Combining \eqref{step1}-\eqref{step2}-\eqref{step3} yields 

\begin{equation}\label{goodpart}
\int_{U} |D \ul-R|^{2}\,dx \leq 
C \int_{U} \dist^{2}(D u,SO(n)) \wedge 
(|D u|^{p}+1) \,dx \,.
\end{equation}

\no
Next we estimate the integral of 
$|D u-R|^{2}\wedge (|D u|^{p}+1) $ in the set $U\setminus U_{\lambda}$. 
In order to 
do this we use again the fact that, for sufficiently large $\lambda$, 
$|D u-R|^{2}\wedge (|D u|^{p}+1)$ is equal to $|D u|^{p}+1$
if $|D u|>\lambda$, and is bounded by a constant if $|D u|<\lambda$.
Hence we write
 
\begin{align}\label{badpart}
\non
\int_{U\setminus U_{\lambda}}\!\! |D u-R|^{2}\!\wedge\! 
(|D u|^{p}+1) \,dx  & \leq 
C |U\setminus U_{\lambda}|+
C \int_{U\setminus U_{\lambda}\cap \{|Du|>\lambda\}}\!\!(|D u|^{p}+1) \,dx  \\
\non
& \leq C\lambda^{-p}\!\int_{\{|D u|>\lambda\}}\!\!|D u|^{p}\,dx \, +
C\!\int_{\{|Du|>\lambda\}}\!\!\!\dist^{2}(D \ul,SO(n))\!\wedge\! 
(|D u|^{p}+1) \,dx    \\
& \leq C\int_{U}\dist^{2}(D \ul,SO(n))
\wedge (|D u|^{p}+1) \,dx \,.
\end{align}
  
\no
Finally \eqref{goodpart} and \eqref{badpart} imply 

\begin{align*}
\int_{U} |D u-R|^{2}\wedge 
(|D u|^{p}+1) \,dx  & =
\int_{U_{\lambda}} |D \ul-R|^{2}\wedge 
(|D \ul|^{p}+1) \,dx +
\int_{U\setminus U_{\lambda}} |D u-R|^{2}\wedge 
(|D u|^{p}+1) \,dx \\
& \leq C\int_{U}\dist^{2}(D \ul,SO(n))
\wedge (|D u|^{p}+1) \,dx \,.
\end{align*}

\end{proof}

The next proposition will be used in the proof of Proposition \ref{propo1}.

\begin{proposition}\label{poincaretype}
Let $n\geq 2$, and let $1\leq p < 2$. Suppose that $U\subset\R^{n}$ 
is a bounded Lipschitz domain. 
Let $u\in W^{1,p}(U,\R^{n})$ satisfy 
$\dsp\int_{U} u \,dx=0$ and  

\begin{equation}\label{small}
\int_{U} |D u|^{2}\wedge (|D u|^{p}+1) \,dx < \e \,,
\end{equation}
\no
with $0<\e<1$. 
Then there exists a constant $C(U,p)$ such that 

\begin{equation}\label{psmall}
\int_{U} (|u|^{2}+|D u|^{2})\wedge (|D u|^{p}+|u|^{p}+1) \,dx 
< C(U,p) \e^{\frac{p}{2}} \,.
\end{equation}
\end{proposition}

\begin{proof}
Let $c>0$ be solution of $c^{2}=c^{p}+1$.
We first provide an estimate for the $L^{p}$ norm of $u$:

\begin{align}\label{smallpnorm}
\int_{U}|D u|^{p} \,dx &\leq 
\int _{\{|D u|\geq c\}}(|D u|^{p} +1) \,dx +
\int _{\{|D u|\leq c\}}
|D u|^{p} \,dx \\
\non &\leq\e+
C\Big(\int_{\{|D u|\leq c\}}|D u|^{2}\,dx
\Big)^{\frac{p}{2}}\\
\non &\leq \e + C\e^{\frac{p}{2}}\\
\non &\leq C\e^{\frac{p}{2}}\,.
\end{align}

\no
Now fix $\lambda>1$ and let $\ul$ be given by Proposition \ref{propolip}. 
Set $U_{\lambda}:=\{u=\ul\}$ and observe that \eqref{smallpnorm} and \eqref{con2} 
imply

\begin{align}\label{badsetest}
\non |U\setminus U_{\lambda}|&\leq \frac{C}{\lambda^{p}}
\int_{\{|D u|>\lambda\}}|D u|^{p}\,dx \\
\non &\leq C\frac{\e^{\frac{p}{2}}}{\lambda^{p}}  \\
&\leq C\e^{\frac{p}{2}}  \,.
\end{align}

\no
Recalling that $\textstyle{\int_{U}u\,dx=0}$, from 
\eqref{smallpnorm}, \eqref{badsetest} and the Poincar\'e inequality 
we deduce that 

\begin{align}
\label{est1}
\int_{U\setminus U_{\lambda}}
(|u|^{2}+|D u|^{2})\wedge
(|D u|^{p}+|u|^{p}+1) \,dx 
\non & \leq
\int_{U\setminus U_{\lambda}}
(|D u|^{p}+|u|^{p}+1) \,dx \\
\non &\leq C \int_{U}|D u|^{p}\,dx + 
|U\setminus U_{\lambda}| \\
&\leq C \e^{\frac{p}{2}}
\end{align}

\no

\no
Next remark that the funtion $\ul-u$ is zero on a set of positive measure.
Therefore the Poincar\'e inequality combined with \eqref{smallpnorm} 
and \eqref{con3}, yields

\begin{align}\label{average}
\Big|-\hspace{-4.5mm}\int_{U}\ul \,dx\Big| & \leq C 
\int_{U\setminus U_{\lambda}}|u-\ul| \,dx \\
\non &\leq C\|D\ul-Du\|_{L^{1}(U)}\\
\non &\leq C\int_{\{|Du|>\lambda\}}|Du| \,dx \\
\non &\leq C \e^{\frac{1}{2}}\,.
\end{align}

\no
Finally, taking into account \eqref{average} and the fact that $u=\ul$ in 
$U_{\lambda}$, we obtain

\begin{align}\label{est2}
\int_{U_{\lambda}}
(|\ul|^{2}+|D\ul|^{2})\wedge(|D\ul|^{p}+|\ul|^{p}+1) \,dx 
\non &\leq 
\int_{U}(|\ul|^{2}+|D\ul|^{2}) \,dx \\
\non &\leq 
C \int_{U}|D\ul|^{2}\,dx + C \e\\
\non &= C \int_{U}|D\ul|^{2}\wedge(|D\ul|^{p}+1)\,dx 
+ C\e \\
\non &\leq C(\e + |U\setminus U_{\lambda}|+ \e) \\
\non &\leq C(\e +\e^{\frac{p}{2}})\\
&\leq C \e^{\frac{p}{2}} \,.
\end{align}

\no
Combining \eqref{est2} with \eqref{est1} yields \eqref{psmall}.

\end{proof}

\begin{lemma}\label{equivalent}
Let $G\in\Mtt$. 
There exist $c_{1},c_{2}>0$ such that for every $A\in\Mtt$

\begin{equation}\label{pointwiseineq}
c_{1}\Big(|A|^{2}\wedge (|A|^{p}+1)\Big) \leq 
|A|^{2}\wedge (|A+G|^{p}+1) \leq
c_{2}\Big(|A|^{2}\wedge (|A|^{p}+1) \Big)\,.
\end{equation} 
 
\end{lemma}

\begin{proof}
We first observe that there exist two positive constants $c_{1},c_{2}$, 
depending on $|G|$ and $p$, such that 
\begin{equation}\label{facile}
c_{1}(|A|^{p}+1) \leq |A+G|^{p}+1\leq
c_{2}(|A|^{p}+1) \,,\quad \forall A\in\Mtt\,.
\end{equation}
Indeed, let us fix $\rho>1$ such that 
$\dsp |A+G|^{p}+1>\frac{1}{2}(|A|^{p}+1)$ for $|A|>\rho$.   
Then, for $\dsp 0\leq|A|\leq\rho$, we have   
$\dsp |A+G|^{p}+1\geq 1\geq \frac{1}{\rho^{p}+1}(|A|^{p}+1)$.
The second inequality of \eqref{facile} is trivial.
In order to prove \eqref{pointwiseineq} it is enough to observe that  
if  $\dsp |A|^{2}\leq |A+G|^{p}+1$ and $\dsp |A|^{p}+1\leq |A|^{2}$, 
then, by \eqref{facile},
$\dsp |A|^{p}+1\leq |A|^{2}\leq c_{2}(|A|^{p}+1)$. 
If otherwise $\dsp|A+G|^{p}+1\leq|A|^{2}$ and $\dsp |A|^{2}\leq |A|^{p}+1$, 
then \eqref{facile} implies
$\dsp c_{1}|A|^{2}\leq |A+G|^{p}+1\leq|A|^{2}$.

\end{proof}

\section{Competition between elastic deformation and formation of dislocations}

\no
We introduce the set
\begin{equation}\label{setc}
\C:=\{F\in L^{p}_{loc}(\R\times S;\Mtt) \::
\curl F = -\b\otimes \dot\Gamma \,d\H^{1}\llcorner \Gamma \}\,.
\end{equation}
The cost associated with a transition of the elastic deformation 
from one well to the other is defined as  

\begin{equation}\label{gammah}
\gh(R, \Gamma):=\inf\Big\{\int_{(-M,M)\times S}W(x,F(x))\,dx \::\:
F\in\C_{M}(R,H), M>0 \Big\}\,,
\end{equation}

\no
where, for each $M>0$, and each $P,Q\in\Mtt$, the set $\C_{M}(P,Q)$ is defined as 

\begin{equation*}
\C_{M}(P,Q):=\{F\in \C\::
F=P \text{ in }(-\infty,-M)\,,\, F=Q \text{ in }(M,+\infty) \}\,.
\end{equation*}

It will be convenient to introduce the quantity $\gh(R,\varnothing)$, defined as the 
minimum cost of a transition in the case when no dislocation is present, i.e., 
$\gh(R,\varnothing)$ is obtained by requiring, in \eqref{setc}, $\curl F=0$ in the 
sense of distributions.

\begin{proposition}\label{invariance} 
For each $\dsp R\in SO(3)$ we have 
$$\gh(R,\varnothing)=\gh(I,\varnothing)\,,\quad
\gh(R,\Gamma)=\gh(I,\Gamma)\,.
$$
\end{proposition}

\no
The proof of Proposition \ref{invariance} can be found in \cite[Proposition 2.4]{mm} for 
the dislocation free case. 
The case with dislocations is treated in a fully analogous way. \\

For ease of notation we set $\gh(\Gamma):=\gh(I,\Gamma)$ and 
$\gh(\varnothing):=\gh(I,\varnothing)$.
Let us remark that such quantities also depend on the radius $r$ of the cross section 
$S$. \\

\no
{\bf Notation.} 
We will write $S_r$, $\gh(\Gamma,r)$ or $\gh(\varnothing,r)$ to emphasize the dependance 
on $r$ when the radius of the cross section plays an essential role.  
In most of the cases however, the dependence on $r$ will be omitted
not to overburden notation. 

\begin{proposition}\label{lowerbound}
Suppose that for each $R\in SO(3)$ and for each $a\in\R^{3}$ 
\begin{equation*}
R-H\neq a\otimes e_{1}\,.
\end{equation*}
Then $\gh(\Gamma)>0$ and $\gh(\varnothing)>0$.
\end{proposition}

\begin{proof}
We will show that $\gh(\Gamma)>0$, the proof for $\gh(\varnothing)$ being completely 
analogous.
By contradiction suppose that $\gh(\Gamma)=0$. 
Then by definition of $\gh(\Gamma)$, there exists 
a sequence $\{F^{(j)}\}\subset\C$ such that 

\begin{equation}\label{infinitesimal}
\int_{\Om}
W(x,F^{(j)}(x))\,dx\to 0\,.
\end{equation}

\no
As already remarked in the introduction, the fields $F^{(j)}$ are locally gradients 
of Sobolev functions. Therefore we can find a set $D\subset S$ and a sequence of functions $\{v^{(j)}\}\subset  W^{1,p}(\Om\setminus D;\R^{3})$ such that 
$\Om\setminus D$ is simply connected and $F^{(j)}=Dv^{(j)}$ in 
$\Om\setminus D$. 
We now apply the rigidity estimate \eqref{prigidity} in combination with  
\eqref{infinitesimal} and the growth condition (iii), to find sequences 
$\{R_{1}^{(j)}\}, \{R_{2}^{(j)}\}\subset SO(3)$ such that 

\begin{eqnarray}\label{goingtozero}
&\dsp \int_{(-L,0)\times S}|Dv^{(j)}-R_{1}^{(j)}|^{2}\wedge 
(|Dv^{(j)}|^{p}+1)\,dx \to 0\,, \\
\non &\dsp\int_{(0,L)\times S}|Dv^{(j)}-R_{2}^{(j)}H|^{2}\wedge 
(|Dv^{(j)}|^{p}+1)\,dx \to 0 \,.
\end{eqnarray}

\no
The first formula of \eqref{goingtozero} implies that 
$|Dv^{(j)}-R_{1}^{(j)}|\to 0$ in measure. Moreover, since  
$\{R_{1}^{(j)}\}$ is a bounded sequence, we have that 
$\{Dv^{(j)}-R_{1}^{(j)}\}$ is bounded in 
$L^{p}\big((-L,0)\times S;\Mtt\big)$, and therefore, up to subsequences (not relabeled),  
$Dv^{(j)}-R_{1}^{(j)}\to 0$ strongly in $L^{q}\big((-L,0)\times S;\Mtt\big)$ for each 
$1\leq q<p$.
Using the second formula of \eqref{goingtozero} and arguing in a similar way for 
the sequence $\{Dv^{(j)}-R_{2}^{(j)}H\}$, we 
deduce that $Dv^{(j)}-R_{2}^{(j)}H\to 0$ strongly in 
$L^{q}\big((0,L)\times S;\Mtt\big)$ for each $1\leq q<p$.
By the Poincar\'e inequality there exist $\dsp\{c_{1}^{(j)}\},\{c_{2}^{(j)}\}\subset\R$ such 
that

\begin{eqnarray*}
& \|v^{(j)}-R_{1}^{(j)}x-c_{1}^{(j)}\|\to 0 
\quad\text{ strongly in }W^{1,q}((-L,0)\times S;\R^{3})\,,\\
& \|v^{(j)}-R_{2}^{(j)}Hx-c_{2}^{(j)}\|\to 0
\quad\text{ strongly in }W^{1,q}((0,L)\times S;\R^{3})\,.
\end{eqnarray*}

\no
Finally, by the trace theorem we find that 
$$
\big(R_{1}^{(j)}-R_{2}^{(j)}H\big)\left(\!\!\!\begin{array}{c}
0\\
x_{2}\\
x_{3}
\end{array}\!\!\!\right)
+c_{1}^{(j)}-c_{2}^{(j)}\to 0  \quad \text{ for }\H^{2}\text{--a.e. } x\in S\setminus D\,,
$$
which yields the contradiction
$\dsp R_{1}-R_{2}H=a\otimes e_{1}$ for some $a\in\R^{3}$, and 
$R_{1},R_{2}\in SO(3)$.

\end{proof}

\begin{remark}\label{scaling}
It can be easily checked that $\gh(\varnothing,r)=r^{3}\gh(\varnothing,1)$. 
Indeed, if $u$ is a competitor for $\gh(\varnothing,1)$, then $u_{r}(x):=ru(x/r)$ 
is a competitor for $\gh(\varnothing,r)$.
\end{remark}

The next proposition provides an upper bound for the energy in the dislocation free case.
We will denote by $\O(\delta)$, with $\delta\in\R$, 
any matrix with norm 
$|\delta|$, i.e.,  $|\O(\delta)|=|\delta|$.

\begin{proposition}\label{deltamatrix}
Assume that $H=I+\O(\delta)$. Then  
$\gh(\varnothing,r)\leq C \delta^{2}r^{3}$.
\end{proposition}

\begin{proof}
Let $u(x):=\var(x_{1})x+(1-\var(x_{1}))Hx$, where 
$\var(x_{1})=1$ for $x_{1}\leq-r/2$, $\var(x_{1})=0$ for $x_{1}\geq r/2$ and 
$\var(x_{1})=-x_{1}/r+1/2$ for $x_{1}\in(-r/2,r/2)$. 
Then one easily checks that 

$$
Du(x)=I+(1-\var(x))\O(\delta)+\frac{1}{r}
\chi_{(-\frac{r}{2},\frac{r}{2})}\O(\delta)x\otimes e_{1}\,.
$$ 
Therefore  $\|Du-I\|_{L^{\infty}}\leq 2|\delta|$, and

\begin{equation*}
\gh(\varnothing,r)\leq
\int_{(-\frac{r}{2},\frac{r}{2})\times S}
W(x,Du(x))\,dx\leq  C\delta^{2}r^{3}\,.
\end{equation*}

\end{proof}

\begin{remark}
Under the assumptions of Proposition \ref{deltamatrix}, we have that 
$\gh(\varnothing,r)\to 0$ as $\delta\to 0$.
\end{remark}
\par

In the next theorem we show that if the radius of the cross section is 
sufficiently large, than formation of dislocations is energetically more convenient than 
purely elastic deformation. 


\begin{theorem}\label{mainthm35}
The following inequality holds

\begin{equation}\label{energeticineq} 
\gh(\varnothing,1)
>
\limsup_{r\to+\infty}\inf_{\Gamma}\frac{\gh(\Gamma,r)}{r^{3}}\,.
\end{equation}

\end{theorem}

\begin{proof}

Let $Q_r\subset\R^2$ be the square of side $2r$ centered at the origin, i.e.,

$$
Q_r:=\{(x_1,x_2)\in\R^2: |x_1|<r,|x_2|<r\}\,.
$$ 

\no
In analogy with $\gh(\varnothing,r)$, we define $\ghh(\varnothing,r)$ 
as the minimum cost associated with 
a transition of the elastic deformation from one well to the other when the cross 
section of the beam is $Q_r$, namely

\begin{align*}
\ghh(\varnothing,r):=
\inf\Big\{&\int_{(-M,M)\times Q_r}W(x,Du(x))\,dx \::\:   \\
& u\in W^{1,p}_{loc}(\R\times Q_r;\R^2)\,,
Du=I \text{ for }x_1<-M, Du=H \text{ for }x_1>M\,, M>0 \Big\}\,.
\end{align*}

\no
We first prove \eqref{energeticineq} for $\ghh$, namely
when the cross-section of the beam is $Q_r$.

From Remark \ref{scaling} it follows that 

\begin{equation}\label{useful}
\gh(\varnothing,r)\leq \ghh (\varnothing,r)\leq \gh(\varnothing,\sqrt{2}r)
= \sqrt{8}\,\gh(\varnothing,r)\,.
\end{equation}

\begin{figure}
\centering
\psfrag{a}{$p_1$}
\psfrag{b}{$p_2$}
\psfrag{c}{$p_3$}
\psfrag{d}{$p_4$}
\psfrag{e}{$\tilde p_1$}
\psfrag{f}{$\tilde p_2$}
\psfrag{g}{$\tilde p_3$}
\psfrag{h}{$\tilde p_4$}
\psfrag{m}{$2\mu$}
\psfrag{r}{$r$}
\includegraphics[scale=0.9]{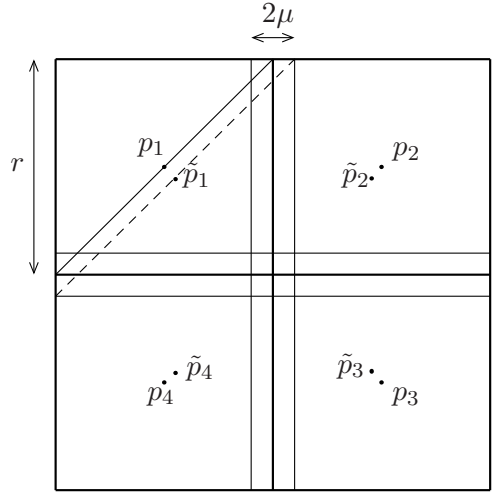}
\caption{The solid and dashed diagonals are of $Q_{\frac{r}{2}}(p_1)$ and 
 $Q_{\frac{r+\mu}{2}}(\tp_1)$ respectively.}
\label{squares}
\end{figure}

\no
We decompose $Q_r$ into the union of four sub-squares of side $r/2$.
Set

$$
p_1=\Big(-\frac{r}{2}\,,\,\frac{r}{2}\Big)\,, \,\,
p_2=\Big(\frac{r}{2}\,,\,\frac{r}{2}\Big)\,, \,\, 
p_3=\Big(\frac{r}{2}\,,\,-\frac{r}{2}\Big)\,, \,\, 
p_4=\Big(-\frac{r}{2}\,,\,-\frac{r}{2}\Big)\,, 
$$

\no
and let $Q_{\frac{r}{2}}(p_1)$, 
$Q_{\frac{r}{2}}(p_2)$, 
$Q_{\frac{r}{2}}(p_3)$, 
$Q_{\frac{r}{2}}(p_4)$ be the 
squares of side $r$ centered at $p_1$, $p_2$, $p_3$ and $p_4$ 
respectively (see Fig. \ref{squares}).
Next we decompose $Q_r$ into the union of four sub-squares overlapping 
on stripes of thickness $2\mu $, with $\mu\ll r$. Specifically, set 

$$
\tp_1=\Big(-\frac{r}{2} + \frac{\mu}{2} \,,\,\frac{r}{2} - \frac{\mu}{2} \Big)\,, \,\,
\tp_2=\Big(\frac{r}{2} - \frac{\mu}{2}\,,\,\frac{r}{2} - \frac{\mu}{2} \Big)\,, \,\, 
\tp_3=\Big(\frac{r}{2} - \frac{\mu}{2} \,,\,-\frac{r}{2} + \frac{\mu}{2} \Big)\,, \,\, 
\tp_4=\Big(-\frac{r}{2} + \frac{\mu}{2}\,,\,-\frac{r}{2} + \frac{\mu}{2}\Big)\,, 
$$

\no
and let $Q_{\frac{r+\mu}{2}}(\tp_1)$, 
$Q_{\frac{r+\mu}{2}}(\tp_2)$, 
$Q_{\frac{r+\mu}{2}}(\tp_3)$, 
$Q_{\frac{r+\mu}{2}}(\tp_4)$ be the 
squares of side $r+\mu$ centered at $\tp_1$, $\tp_2$, $\tp_3$ and 
$\tp_4$ respectively (see Fig. \ref{squares}).
Now fix $\delta>0$. 
By definition of $\ghh(\varnothing,\frac{r+\mu}{2})$, 
there exist $M>0$ and 
$u\in W^{1,p}_{loc}\big(\R\times Q_{\frac{r+\mu}{2}};\R^3\big)$ such 
that  $u=x$  in $(-\infty,-M)\times Q_{\frac{r+\mu}{2}}$,  
$u=Hx$ in $(M,+\infty)\times Q_{\frac{r+\mu}{2}}$ and 

\begin{equation}\label{eq:basic}
\int_{(-M,M)\times Q_{\frac{r+\mu}{2}}}W(x,Du(x))\,dx = 
\ghh\big(\varnothing,\textstyle\frac{r+\mu}{2}\big) + \delta \,.
\end{equation}

\no
Up to an arbitrarily small error in \eqref{eq:basic}, 
by applying Proposition \ref{propolip} we can assume that 
$u\in W^{1,\infty}\big((-M,M)\times Q_{\frac{r+\mu}{2}};\R^3\big)$.
(Remark that $M=\O(r)$.)
Define

\begin{align*}\dsp
& u_1(x):=\,\,  u\big(x -(0,\tp_1)\big) \hspace{3.2cm}
\text{ in } \, (-\infty,+\infty)\times Q_{\frac{r+\mu}{2}}(\tp_1)\,,\\
& u_2(x):= \begin{cases}
u\big(x -(0,\tp_2)\big) + (0,\tp_2-\tp_1) &  
                                    \text{ in }(-\infty,0)\times Q_{\frac{r+\mu}{2}}(\tp_2)\,,\\
u\big(x -(0,\tp_2)\big) + H(0,\tp_2-\tp_1) &  
                                     \text{ in }(0,+\infty)\times Q_{\frac{r+\mu}{2}}(\tp_2)\,,
\end{cases}   \\                                  
& u_3(x):= \begin{cases}
u\big(x -(0,\tp_3)\big) + (0,\tp_3-\tp_1)  & \text{ in }
                                                      (-\infty,0)\times Q_{\frac{r+\mu}{2}}(\tp_3)\,,\\
u\big(x -(0,\tp_3)\big) + H(0,\tp_3-\tp_1)   &  
                                     \text{ in }
                                     (0,+\infty)\times Q_{\frac{r+\mu}{2}}(\tp_3)\,,
\end{cases}   \\     
& u_4(x):= \begin{cases}                                 
u\big(x -(0,\tp_4)\big) + (0,\tp_4-\tp_1)  &  \text{ in }
                                                      (-\infty,0)\times Q_{\frac{r+\mu}{2}}(\tp_4)\,,\\
u\big(x -(0,\tp_4)\big) + H(0,\tp_4-\tp_1)   &  
                                     \text{ in }
                                     (0,+\infty)\times Q_{\frac{r+\mu}{2}}(\tp_4)\,,
\end{cases}
\end{align*}

\no
where, abusing notation, 
$(0,\tp_i)$ denotes the point $(0,(\tp_i)_1,(\tp_i)_2)$, for $i=1,\dots,4$.
The constants in the definition of $u_i$ have been chosen so as to 
ensure that $D u_i=I$ for $x_1<-M$ and $D u_i=H$ for $x_1>M$ for all 
$i=1,\dots,4$. The function $u_i$, for $i=2,3,4$, has a constant jump equal to 
$(H-I)(0,\tp_i-\tp_1)$ on the set $\{0\}\times Q_{\frac{r+\mu}{2}}(\tp_i)$.
\no
Introduce the function 

\begin{equation*}
v(x):= \begin{cases}
u_1(x) & \text{ for } x \in (-\infty,+\infty)\times Q_{\frac{r}{2}}(p_1)\,,\\
u_2(x) & \text{ for } x \in (-\infty,+\infty)\times Q_{\frac{r}{2}}(p_2)\,,\\
u_3(x) & \text{ for } x \in (-\infty,+\infty)\times Q_{\frac{r}{2}}(p_3)\,,\\ 
u_4(x) & \text{ for } x \in (-\infty,+\infty)\times Q_{\frac{r}{2}}(p_4)\,.
\end{cases}                                  
\end{equation*}

\no
We modify $v$ in order to remove the jump on the sets $\{x_2 =0\}$ 
and $\{x_3=0\}$. 
Introduce cylindrical coordinates 

$$
x_1=\rho\cos\theta\,, \quad  x_2=\rho\sin\theta\,,
$$

\no
and define the sector 

$$
\omega_1:=\Big\{(x_1,x_2,x_3)\in\R\times Q_r: x_1=\rho\cos\theta\,, x_2=\rho\sin\theta\,, |\tan\theta|< \frac{\mu}{M}\,,
\rho\in(-M,M) \Big\}\,.
$$

\begin{figure}
\centering
\psfrag{m}{$2\mu$}
\psfrag{o}{$\omega_1$}
\includegraphics{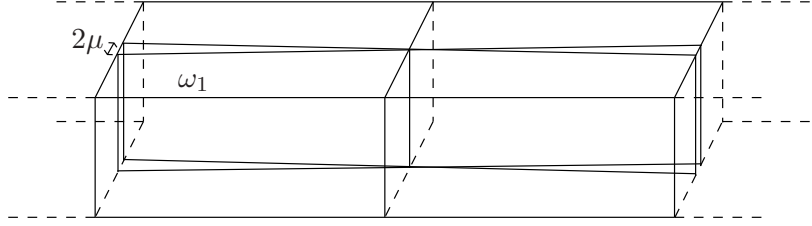}
\caption{The set $\omega_1$.}
\label{set-omega}
\end{figure}

\no
Let $\dsp\bar\theta=\arctan\frac{\mu}{M} $, and define

\begin{equation*}
v^*(x):= \begin{cases} 
\vspace{2mm}

u_1(x) + \frac{\tan\bar\theta -\frac{x_2}{x_1}}{2\tan\bar\theta}
\big(u_2(x)-u_1(x)\big) 
& \text{ for } x \in \omega_1\,, x_1 < 0  \,, x_3 > 0 \\

\vspace{2mm}

u_4(x) + \frac{\tan\bar\theta -\frac{x_2}{x_1}}{2\tan\bar\theta}
\big(u_3(x)-u_4(x)\big) 
& \text{ for } x \in \omega_1\,, x_1 < 0 \,, x_3 < 0 \\

\vspace{2mm}

u_2(x) + \frac{\tan\bar\theta -\frac{x_2}{x_1}}{2\tan\bar\theta}
\big(u_1(x)-u_2(x)\big) 
& \text{ for } x \in \omega_1\,, x_1 > 0  \,, x_3 > 0 \\

\vspace{2mm}

u_3(x) + \frac{\tan\bar\theta -\frac{x_2}{x_1}}{2\tan\bar\theta}
\big(u_4(x)-u_3(x)\big) 
& \text{ for } x \in \omega_1\,, x_1 > 0 \,, x_3 < 0 \\

\vspace{2mm}
v(x) & \text{ for } x \in \R\times Q_r\setminus\omega_1\,.
\end{cases} 
\end{equation*}

\no
By an analogous interpolation we further modify  $v^*$ in the 
sector 

$$
\omega_2:=\Big\{(x_1,x_2,x_3)\in\R\times Q_r: x_1=\rho\cos\theta\,, x_3=\rho\sin\theta\,, |\tan\theta|< \frac{\mu}{M}\,,
\rho\in(-M,M) \Big\}\,,
$$

\no
thus obtaing a function 
$\dsp v^{**} \in 
W^{1,p}_{loc}
\Big(
\big(\R\times Q_r\big)\setminus \big(\{0\}\times 
\cup_{i=2}^4 Q_{\frac{r}{2}}(p_i)\big) ;\R^3
\Big)$, 
such that $v^{**}$ has a constant jump on the set 
$\{0\}\times Q_{\frac{r}{2}}(p_i)$,  for $i=2,3, 4$. 
(Remark that the condition $Dv^{**}\in L^p$ is ensured by the boundedness 
of $u$ in the sets $\omega_1$ and $\omega_2$.)
The dislocation line $\Gamma(v^{**})$ associated with $v^{**}$ is the boundary of the jump set of $v^{**}$, i.e., 

$$
\Gamma(v^{**})=\large\bigcup_{i=2}^4\partial Q_{\frac{r}{2}}(p_i)  \cap Q_r \,.
$$

\no
Notice that each curve $\partial Q_{\frac{r}{2}}(p_i)  \cap Q_r $ is associated 
with a different Burgers vector $b_i$, specifically $b_i=(H-I)(0,\tp_i-\tp_1)$. 
Using \eqref{eq:basic}, one finds 

\begin{equation}\label{eq:square-estimate}
\begin{aligned}
\inf_{\Gamma}\ghh(\Gamma,r) \leq \ghh(\Gamma(v^{**}),r)
& \leq
\int_{(-M,0)\times Q_r} W(x,Dv^{**})\, dx 
+
\int_{(0,M)\times Q_r} W(x,Dv^{**})\, dx \\
& \leq 
4 \Big(\ghh\big(\varnothing,\textstyle\frac{r+\mu}{2}\big) + \delta\Big) +
\dsp \int_{\omega_1\cup\omega_2} W(x,Dv^{**})\, dx  \,.
\end{aligned}
\end{equation}

\no
By a suitable choice of $\mu=o(r)$, one can prove that 
$\dsp    \int_{\omega_1\cup\omega_2} W(x,Dv^{**})\, dx = o(r^3)$. Then, 
since $\delta$ is arbitrary,  from 
\eqref{eq:square-estimate} we deduce 

$$
\limsup_{r\to+\infty}\inf_{\Gamma}\frac{\gh(\Gamma,r)}{r^{3}}
\leq
\frac{1}{2}\ghh(\varnothing,1)\,.
$$

\no
This proves  \eqref{energeticineq} for $\ghh$, 
but it is not enough to prove it for $\gh$. 
Indeed, using \eqref{useful} and taking the 
restriction of $v^{**}$ to $(-\infty,+\infty)\times S_{r}$ we would get 

\begin{equation*}
\begin{aligned}
\gh(\Gamma(v^{**}),r)& \leq \ghh(\Gamma(v^{**}),r) \leq
\frac{r^3}{2}\ghh(\varnothing,1) + o(r^3) + \delta \\
& \leq
\frac{r^3}{2}\gh(\varnothing,\sqrt{2}) + o(r^3) + \delta  \leq
\frac{r^3}{2}\sqrt{8} \ \gh(\varnothing,1) + o(r^3) + \delta \\
& = \sqrt{2}r^3\gh(\varnothing,1) + o(r^3) + \delta  \,.
\end{aligned}
\end{equation*}

\no
In order to decrease the coefficient in front of $\gh(\varnothing,1)$, it is enough 
to divide $Q_r$ into eight sub-squares of side $r/4$ 
(instead of four sub-squares of side $r/2$ as we did before)  
and to repeat the same construction as before, namely, to glue together  suitable 
translations 
of the function $u$, where $u$ is such that 

$$
\int_{(-M,M)\times Q_{\frac{r+\mu}{4}}}W(x,Du(x))\,dx = 
\ghh\big(\varnothing,\textstyle\frac{r+\mu}{4}\big) + \delta \,.
$$

\no
This yields 

$$
\limsup_{r\to+\infty}\inf_{\Gamma}\frac{\gh(\Gamma,r)}{r^{3}}
\leq \frac{\sqrt{2}}{8}\gh(\varnothing,1)\,.
$$

\end{proof}

\section{Compactness and lower bound}\label{section:four}

\begin{theorem}\label{thm1}
Assume $W$ satisfies {\rm (i)-(iii)}. Let $\{\Fh\}\subset \C$ be a 
sequence such that 

\begin{equation}\label{equibounded}
\limsup_{h\to 0^{+}}\frac{1}{h}\int_{\Om}
W(x,\Fh_{h}(x))\,dx\leq c \,.
\end{equation}

\no
Then there exists a subsequence (not relabeled) such that 

\begin{equation*}
\Fh_{h}\weak F \quad \text{ weakly in } L^{p}(\Om,\Mtt)\,,
\end{equation*}

\no
where $F$ is independent of $x_{2}$ and $x_{3}$ and satisfies

\begin{equation}\label{charact}
F\in 
\begin{cases}
{\rm co}(SO(3))     &  \text{a.e in }(-L,0)\,,\\
{\rm co}(SO(3)H)   &  \text{a.e in }(0,L)\,. 
\end{cases}
\end{equation}

\no
Moreover, for each such subsequence we have

\begin{equation}\label{gammaliminf}
\liminf_{h\to 0^{+}}\frac{1}{h}\int_{\Om}W(x,\Fh_{h}(x))\,dx
\geq\gh \,.
\end{equation}
\end{theorem}

\begin{proof}
\no
The assumption \eqref{equibounded} together with the growth condition from 
below on $W$ imply that the sequence $\{\Fh_{h}\}$ is uniformly bounded in 
$L^{p}(\Om,\Mtt)$. Therefore there exists a subsequence (not relabeled) 
converging to some $F$ weakly in $L^{p}(\Om,\Mtt)$.

In order to prove \eqref{charact} we use the gradient structure of 
$\Fh$ in 
simply connected domains. 
Indeed, as already remarked in the introduction, there exists 
$\vh\in W^{1,p}(\Om\setminus D,\R^{3})$ such that 
$\Fh=D\vh$ a.e. in $\Om\setminus D$.
Next we divide the intervals $(-L,0)$ and $(0,L)$ into subintervals of 
length $\tau_{h}\sim h$ and apply the rigidity estimate \eqref{prigidity} 
to $\uh(z):=\vh(z_{1},\frac{z_{2}}{h},\frac{z_{3}}{h})$ in the Cartesian product of each subinterval 
and the cross-section $hS$. More precisely, let $\tau_{h}:=L/[L/h]$ where $[t]$ denotes 
the largest 
integer less than or equal to $t$. 
Then, for every $h>0$ and $a\in [-L,0)\cap \tau_{h}\Z$ there exists $G^{h}(a)\in SO(3)$ such that  

\begin{equation*}
\int_{(a,a+\tau_{h})\times hS} 
|D\uh-G^{h}(a)|^{2}\wedge (|D\uh|^{p}+1) \,dz \leq C
\int_{(a,a+\tau_{h})\times hS} 
\dist^{2}(D\uh,SO(3))\wedge (|D\uh|^{p}+1) \,dz \,,
\end{equation*}

\no
and, for every $h>0$ and $a\in [0,L)\cap \tau_{h}\Z$, there exists $G^{h}(a)\in SO(3)H$ such that

\begin{equation*}
\int_{(a,a+\tau_{h})\times hS} 
|D\uh-G^{h}(a)|^{2}\wedge (|D\uh|^{p}+1) \,dz \leq C
\int_{(a,a+\tau_{h})\times hS} 
\dist^{2}(D\uh,SO(3)H)\wedge (|D\uh|^{p}+1) \,dz \,.
\end{equation*}

\no
By interpolation one defines a piecewise constant matrix field 
$G^{h}:(-L,L)\to SO(3)\cup SO(3)H$ such that 
$G^{h}(x_{1})=G^{h}(a)$ if $x_{1}\in (a,a+\tau_{h})$ and 
$a\in (-L,L)\cap \tau_{h}\Z$. By rescaling the problem back to $\Om$, one gets

\begin{equation*}
\int_{(-L,0)\times S}\!\! 
|\Fh_{h}(x)-G^{h}(x_{1})|^{2}\wedge (|\Fh_{h}(x)|^{p}+1) \,dx \leq C\!\!
\int_{(-L,0)\times S} \!\!
\dist^{2}(\Fh_{h}(x),SO(3))\wedge (|\Fh_{h}(x)|^{p}+1) \,dx \,,
\end{equation*}

\begin{equation*}
\int_{(0,L)\times S}\!\! 
|\Fh_{h}(x)-G^{h}(x_{1})|^{2}\wedge (|\Fh_{h}(x)|^{p}+1) \,dx \leq C \!\!\\
\int_{(0,L)\times S} \!\!
\dist^{2}(\Fh_{h}(x),SO(3)H)\wedge (|\Fh_{h}(x)|^{p}+1) \,dx \,.
\end{equation*}

\no
The above inequalities and \eqref{equibounded} imply that 
$|\Fh_{h}-G^{h}|\to 0$ in measure 
and therefore $\Fh_{h}-G^{h}\to 0$ in $L^{q}(\Om;\Mtt)$ for each $q<p$, 
the sequence $\{\Fh_{h}-G^{h}\}$ being uniformly bounded in $L^{p}(\Om;\Mtt)$.
On the other hand, the sequence $\{G^{h}\}$ is uniformly bounded in 
$L^{\infty}((-L,L);\Mtt)$ and therefore, up to subsequences, it converges
to some $G\in L^{\infty}((-L,L);\Mtt)$ in the weak$^{*}$ topology of 
$L^{\infty}((-L,L);\Mtt)$. 
Then \eqref{charact} easily follows from the fact that $G^{h}\in SO(3)$ a.e. 
in $(-L,0)$ and $G^{h}\in SO(3)H$ a.e. in $(0,L)$.

In order to show \eqref{gammaliminf}, we define 

$$
\tilde{F}^{(h)}(x):=\Fh_{h}(hx_{1},x_{2},x_{3})\,,
$$

\no
and observe that $\tilde{F}^{(h)}\subset\C$, since 
$\dsp\tilde{F}^{(h)}(x)=D\Big(\frac{1}{h}\uh(hz)\Big)$. 
Moreover 

\begin{equation*}
\frac{1}{h}\int_{\Om}W(x,\Fh_{h}(x))\,dx  =
\int_{(-L_{h},L_{h})\times S}W(x,\tilde{F}^{(h)}(x))\,dx \,,
\end{equation*}

\no
where $\dsp L_{h}:=\frac{L}{h}$.
Then by \eqref{equibounded} we have 

\begin{equation*}
\int_{(-L_{h},0)\times S} 
\dist^{2}(\tilde{F}^{(h)},SO(3))\wedge (|\tilde{F}^{(h)}|^{p}+1) \,dx \,+\,
\int_{(0,L_{h})\times S} 
\dist^{2}(\tilde{F}^{(h)},SO(3)H)\wedge (|\tilde{F}^{(h)}|^{p}+1) \,dx \leq C \,.
\end{equation*}

\no
Finally Proposition \ref{propo1} below yields

\begin{equation*}
\int_{(-L_{h},L_{h})\times S}W(x,\tilde{F}^{(h)}(x))\,dx
\geq \gh - C \Big(\frac{h}{L}\Big)^{\frac{p}{2}}\,.
\end{equation*}

\end{proof}

\begin{proposition}\label{propo1}
There exists $C>0$ such that for all $M\geq 2$ and all $F\in \C$ 
the following implication holds: if

\begin{equation*}
\int_{(-M,-1)\times S} 
\dist^{2}(F,SO(3))\wedge (|F|^{p}+1) \,dx +
\int_{(1,M)\times S} 
\dist^{2}(F,SO(3)H)\wedge (|F|^{p}+1) \,dx \leq C_{0} \,,
\end{equation*}

\no
then
\begin{equation*}
\int_{(-M,M)\times S} W(x,F) \,dx \geq 
\gh-C\Big(\frac{C_{0}}{M}\Big)^{\frac{p}{2}} \,.
\end{equation*}
\end{proposition}

\begin{proof}
There exists $j\in\{1,\dots,M-1\}$ such that

\begin{equation*}
\int_{(j,j+1)\times S} 
\dist^{2}(F,SO(3)H)\wedge (|F|^{p}+1) \,dx  \leq \frac{C_{0}}{M} \,.
\end{equation*}

\no
Proposition \ref{rigext} implies that there exists $R\in SO(3)$ such 
that 

\begin{equation*}
\int_{(j,j+1)\times S} 
|F-RH|^{2}\wedge (|F|^{p}+1) \,dx \leq \frac{C_{0}}{M} \,,
\end{equation*}

\no
and therefore, by Lemma \ref{equivalent}, we deduce that 

\begin{equation}\label{estimate1}
\int_{(j,j+1)\times S} 
|F-RH|^{2}\wedge (|F-RH|^{p}+1) \,dx \leq \frac{C_{0}}{M} \,.
\end{equation}

\vspace{1mm}

\no 
Let $v\in W_{loc}^{1,p}((\R\times S)\setminus D,\R^{3})$ be such that 
$F=Dv$ a.e. in $(\R\times S)\setminus D$. 
Applying Proposition \ref{poincaretype} to the function $u=v-(RHx+c)$, 
for a suitable $c\in\R^{3}$, we deduce from 
\eqref{estimate1}

\begin{equation}\label{estimate2}
\int_{(j,j+1)\times S}
(|v-(RHx+c)|^{2}+|D v-RH|^{2})
\wedge (|v-(RHx+c)|^{p}+|D v-RH|^{p}+1) 
\,dx\leq C\, \Big(\frac{C_{0}}{M}\Big)^{\frac{p}{2}}.
\end{equation}

\no
Now let 
$\dsp\var\in C^{\infty}(\R)$ be a cut-off function, i.e., $\var=1$ in 
$(-\infty,0)$, $\var=0$ in $(1,+\infty)$, $0\leq\var\leq 1$, and set

$$
\tv(x):=\var(x_{1}-j)v(x)+(1-\var(x_{1}-j))(RHx+c)\,.
$$ 
It is readily seen that $\tv=v$ in $(0,j)\times S$, $\tv=RHx+c$ in 
$(j+1,+\infty)\times S$, and

\begin{equation}\label{gradbound}
|D \tv-RH|\leq C(|v-(RHx+c)|+|D v-RH|)\quad \text{ for }
x_{1}\in(j,j+1) \,.
\end{equation}

\no
Taking into account \eqref{gradbound} and the upper growth condition (iii) 
on $W$, we find

\begin{align}\label{energyest}
\non\int_{(j,j+1)\times S}W(D \tv)\,dx &\leq
\int_{(j,j+1)\times S}
\dist^{2}(D \tv,SO(3)H)\wedge (|D \tv|^{p}+1) \,dx\\
\non & \leq \int_{(j,j+1)\times S}
|D \tv-RH|^{2}\wedge (|D \tv|^{p}+1) \,dx\\
\non &\leq C\int_{(j,j+1)\times S}
(|v-(RHx+c)|^{2}+|D v-RH|^{2})
\wedge \\
&\hspace{2.55cm}(|v-(RHx+c)|^{p}+|D v-RH|^{p}+1) 
\,dx\,.
\end{align}

\no
Combining \eqref{estimate2} and \eqref{energyest} yields

\begin{equation}\label{comparison}
\int_{(0,+\infty)\times S}W(D \tv)\,dx \leq
\int_{(0,j)\times S}W(D v)\,dx + 
C\Big(\frac{C_{0}}{M}\Big)^{\frac{p}{2}}\,.
\end{equation}

\vspace{2mm}

\no
Modifying $v$ in a similar way in some 
subset $(-j'-1,-j')\subset(-M,-1)\times S$, $j'\in\{1,\dots,M-1\}$, 
one defines a map $\tv$ such that 

\begin{equation*}
\tv(x)= 
\begin{cases}
R'x+c'   &  \text{ in }(-\infty,-j'-1)\times S\,,\\
v(x)      &  \text{ in }(-j',j)\times S\,,\\
RHx+c    &  \text{ in }(j+1,+\infty)\times S\,. 
\end{cases}
\end{equation*}

\no
Next set $\tilde{F}:=\nabla \tv$, i.e., $\tilde{F}$ is defined as the absolutely 
continuous part of the gradient $D\tv$, and remark that $\tilde{F}=F$ in 
$(-j',j)\times S$ and $\tilde{F}\in\C_{M}(R',RH)$.
From \eqref{comparison} and the definition of $\gh$, it follows that 

\begin{equation*}
\gh\leq
\int_{(-M,M)\times S}W(\tilde{F})\,dx \leq
\int_{(-M,M)\times S}W(F)\,dx + 
C\Big(\frac{C_{0}}{M}\Big)^{\frac{p}{2}}\,.
\end{equation*}

\end{proof}

\section{upper bound}\label{section:five}

\no
We will use the notation $\dsp x'=\left(
\begin{array}{cc}x_{2}\\ x_{3}
\end{array}
\right)$.

\begin{theorem}\label{thm2}
Assume $W$ satisfies {\rm (i)-(iii)}. 
Let $F\in L^{p}((-L,L);\Mtt)$ satisfy 
 
\begin{equation}\label{realizable}
F\in 
\begin{cases}
{\rm co}(SO(3))     &  \text{a.e in }(-L,0)\,,\\
{\rm co}(SO(3)H)   &  \text{a.e in }(0,L)\,. 
\end{cases}
\end{equation} 

\no
Then there exists a sequence $\{\Fh\}\subset\F^{(h)}$ such that 

\begin{equation}\label{recovery}
\Fh_{h}\weak F \quad\text{ weakly in } L^{p}(\Om;\Mtt)\,,
\end{equation}

\no
and
\begin{equation}\label{gammalimsup}
\limsup_{h\to 0^{+}}\frac{1}{h}\int_{\Om}W(x,\Fh(x))\,dx
\leq\gh \,.
\end{equation}
\end{theorem}

\begin{proof}
The proof is very similar to that in \cite{mm}, therefore we  
refer to \cite[Theorem 3.1]{mm} for full details. 
We first assume that $F$ is piecewise constant with values 
in $K$, i.e., $F\in SO(3)$ for a.e. $x_{1}\in (-L,0)$, and 
$F\in SO(3)H$ for a.e. $x_{1}\in (0,L)$.
In this case there exist $a_{0}=-L<a_{1}<\dots<a_{n+1}=0$, 
$0=b_{0}<\dots<b_{k+1}=L$ 
and $R_{i}\in SO(3)$, $S_{j}\in SO(3)H$, $i=0,\dots,n$, 
$j=0,\dots,k$, such that 

\begin{equation}\label{piececonst}
F=\sum_{i=0}^{n}\chi_{(a_{i},a_{i+1})}R_{i}+
 \sum_{j=0}^{k}\chi_{(b_{j},b_{j+1})}S_{j}\,.
\end{equation}

\no
Let $\{\sigma_{h}\}$ be a sequence of positive numbers such that 
$h\ll\sigma_{h}\ll 1$. We define a function $y^{(h)}$ in the 
following way:

\begin{equation*}
y^{(h)}(x):= 
\begin{cases}
R_{0}\left(\!\!\!\begin{array}{c}x_{1}\\ hx' \end{array}\!\!\!\right)
&  \text{ if } x\in (-L,a_{1}-\sigma_{h})\times S\,,\\
R_{i}\left(\!\!\!\begin{array}{c}x_{1}\\ hx' \end{array}\!\!\!\right) +c_{i}^{(h)}   
&  \text{ if } x\in (a_{i}+\sigma_{h},a_{i+1}-\sigma_{h})\times S\,,
\quad i=1,\dots,n-1\\
R_{n}\left(\!\!\!\begin{array}{c}x_{1}\\ hx' \end{array}\!\!\!\right) +c_{n}^{(h)}    
&  \text{ if } x\in (a_{n}+\sigma_{h},-\sigma_{h})\times S\,, \\
S_{0}\!\left(\!\!\!\begin{array}{c}x_{1}\\ hx' \end{array}\!\!\!\right)
&  \text{ if } x\in (\sigma_{h},b_{1}-\sigma_{h})\times S\,,\\
S_{j}\!\left(\!\!\!\begin{array}{c}x_{1}\\ hx' \end{array}\!\!\!\right) +d_{j}^{(h)}   
& \text{ if } x\in (b_{j}+\sigma_{h},b_{j+1}-\sigma_{h})\times S\,,
\quad j=1,\dots,k-1\,,\\
S_{k}\!\left(\!\!\!\begin{array}{c}x_{1}\\ hx' \end{array}\!\!\!\right) +d_{k}^{(h)}    
&  \text{ if } x\in (b_{k}+\sigma_{h},L)\times S\,,
\end{cases}
\end{equation*} 

\no
where the constants $c_{i}$, $i=1,\dots,n$, and $d_{j}$, 
$j=1,\dots,k$, will be chosen later. 

In order to define $y^{(h)}$ in the set 
$(-\sigma_{h},\sigma_{h})\times S$, we proceed in the 
following way. 
Let $\eta>0$. By definition of $\gh$, \eqref{gammah},  
there exist $M>0$ and $F\in\C_{M}(R_{n},S_{0})$ such that 

\begin{equation*}
F=R_{n} \quad\text{ a.e. in } (-\infty,-M)\,, \quad 
F=S_{0} \quad\text{ a.e. in } (M,+\infty)\,, 
\end{equation*}

\no
and

\begin{equation}\label{goodtest}
\int_{(-M,M)\times S}W(x,F(x))\,dx \leq \gh + \eta \,.
\end{equation}

\no
Now let $v\in W_{loc}^{1,p}((\R\times S)\setminus D,\R^{3})$ be such that 
$F=Dv$ a.e. in $(\R\times S)\setminus D$ and define $y^{(h)}$ in the set 
$((-\sigma_{h},\sigma_{h})\times S)\setminus D$ as 

\begin{equation*}
y^{(h)}(x):= hv\Big(\frac{x_{1}}{h},x'\Big) + l_{0}^{(h)}\,,
\end{equation*}

\no
where the constant $l_{0}^{(h)}$ will be chosen later.

In the sets $(a_{i}-\sigma_{h},a_{i}+\sigma_{h})\times S$, for 
$i=1,\dots,n$, we define $y^{(h)}$ in the following way. 
We construct a smooth function 
$P_{i}:\R\to SO(3)$ such that $P_{i}(0)=R_{i-1}$ and 
$P_{i}(1)=R_{i}$ and we set 
$P_{i}^{(h)}(x_{1}):=P_{i}\Big(
\frac{x_{1}-a_{i}+\sigma_{h}}{2\sigma_{h}}\Big)$.
Then we define for 
$x\in(a_{i}-\sigma_{h},a_{i}+\sigma_{h})\times S$

\begin{equation}\label{interpolation}
y^{(h)}(x):=\int_{a_{i}-\sigma_{h}}^{x_{1}}
P^{(h)}_{i}(s) e_{1} \,ds +
P_{i}^{(h)}(x_{1})
\left(\!\!
\begin{array}{c}0 \\ hx'
\end{array}
\!\!\right)
+l_{i}^{(h)}\,,
\end{equation}

\no
where the constants $l_{i}^{(h)}$ will be chosen later.

In the sets $(b_{j}-\sigma_{h},b_{j}+\sigma_{h})\times S$, for 
$j=1,\dots,k$, we can construct a smooth function 
$Q_{j}:\R\to SO(3)H$ such that $Q_{j}(0)=S_{j-1}$ and 
$Q_{j}(1)=S_{j}$. Then one defines $y^{(h)}$ in 
$(b_{j}-\sigma_{h},b_{j}+\sigma_{h})\times S$ 
as in \eqref{interpolation}.

Next we choose the constants $c_{i},d_{i},l_{i}$ so that, for $h$ sufficiently small, 
the function $y^{(h)}$ belongs to 
$W^{1,p}((\R\times S)\setminus D,\R^{3})$. 
Finally we set $\Fh(x):=Dy^{(h)}(x)$, for a.e. $x\in\Om$. 
As far as \eqref{recovery} is concerned, 
it can be easily checked that in fact $\Fh_{h}$ converges strongly 
to $F$ in $L^{p}(\Om,\R^{3})$.
This follows from the fact that $F^{(h)}$ coincides with $F$ outside 
the sets
$(a_{i}-\sigma_{h},a_{i}+\sigma_{h})\times S$, 
$i=1\dots,n+1$, and 
$(b_{j}-\sigma_{h},b_{j}+\sigma_{h})\times S$, 
$j=1\dots,k$, and is uniformly bounded in such sets 
(see \cite{mm} for further details).

Finally we prove \eqref{gammalimsup}. 
Since $W(x,\Fh_{h})=0$ in the complement of the sets 
$(a_{i}-\sigma_{h},a_{i}+\sigma_{h})\times S$, 
$i=1\dots,n+1$, and 
$(b_{j}-\sigma_{h},b_{j}+\sigma_{h})\times S$, 
$j=1\dots,k$, we have that 

\begin{align*}
\frac{1}{h}\int_{\Om}W(x,\Fh_{h})\,dx = & 
\frac{1}{h}
\int_{(-\sigma_{h},\sigma_{h})\times S}
\!\! W(x,\Fh_{h})\,dx \:+ \\
&\frac{1}{h}\sum_{i=1}^{n}
\int_{(a_{i}-\sigma_{h},a_{i}+\sigma_{h})\times S}
\!\! W(x,\Fh_{h})\,dx +
\frac{1}{h}\sum_{j=1}^{k}
\int_{(b_{j}-\sigma_{h},b_{j}+\sigma_{h})\times S}
\!\! W(x,\Fh_{h})\,dx \,.
\end{align*}

\no
From \eqref{goodtest} it follows that 

\begin{equation}\label{smallenergy1}
\frac{1}{h}
\int_{(-\sigma_{h},\sigma_{h})\times S}
\!\! W(x,\Fh_{h})\,dx \leq \gh +\eta \,.
\end{equation}

\no
On the other hand, using \eqref{interpolation} and 
the growth conditions from above on $W$, 
for each $i=1\dots,n$, we find

\begin{align}\label{smallenergy2}
\non
\frac{1}{h}
\int_{(a_{i}-\sigma_{h},a_{i}+\sigma_{h})\times S}
\!\! W(x,\Fh_{h})\,dx & \leq 
\frac{C_{2}}{h} 
\int_{(a_{i}-\sigma_{h},a_{i}+\sigma_{h})\times S}
\Big(\dist^{2}(\Fh_{h},SO(3))\Big)\wedge 
\Big(|\Fh_{h}|^{p}+1\Big)\,dx \\
\non & \leq 
\frac{C_{2}}{h} 
\int_{(a_{i}-\sigma_{h},a_{i}+\sigma_{h})\times S}
\dist^{2}(\Fh_{h},SO(3))\,dx\\
\non & \leq 
\frac{C_{2}}{h} 
\int_{(a_{i}-\sigma_{h},a_{i}+\sigma_{h})\times S}
|\Fh_{h} - P_{i}^{(h)}|^{2} \,dx \\
& \leq C h \sigma_{h}^{-1} \,.
\end{align}

\no
A similar estimate holds in the set 
$(b_{j}-\sigma_{h},b_{j}+\sigma_{h})\times S$ for each 
$j=1,\dots,k$\,,

\begin{equation}\label{smallenergy3}
\frac{1}{h}
\int_{(b_{j}-\sigma_{h},b_{j}+\sigma_{h})\times S}
\!\! W(x,\Fh_{h})\,dx \leq C h \sigma_{h}^{-1} \,.
\end{equation}

\no
Summing up together 
\eqref{smallenergy1}-\eqref{smallenergy2}-\eqref{smallenergy3}, 
and recalling that $h\sigma^{-1}_{h}\to 0$, we conclude that

\begin{equation*}
\limsup_{h\to 0^{+}}
\frac{1}{h}
\int_{\Om}
\!\! W(x,\Fh_{h})\,dx =\gh +\eta \,.
\end{equation*}

\no
For the general case when \eqref{realizable} holds, one finds 
a sequence of piecewise constant maps $\{F_{j}\}$ as in 
\eqref{piececonst} such that $F_{j}\weak F$ weakly in 
$L^{p}(\Om;\Mtt)$, and then argues by approximation 
(see \cite{mm} for further details).

\end{proof}

It will be convenient to define the following set

\begin{multline*}
\F:=\big\{
F\in L^{p}(\Om,\Mtt):\\
\,F(x)=F(x_{1}),
|F^{1}|\leq 1 \text{ a.e. in }(-L,0), 
|F^{1}|\leq\zeta_{1}\text{ a.e. in }(-L,0),
F^{2}=F^{3}=0
\big\}\,.
\end{multline*}

\begin{remark}
If $F\in\F$, then there exists a map $u\in W^{1,\infty}((-L,L);\R^{3})$, 
such that $F^{1}=u'(x_{1})$.
\end{remark}

The next theorem states that the domain of the $\Gamma$-limit 
of the sequence $\{\frac{1}{h}\I^{(h)}\}$ is $\F$ and that the 
$\Gamma$-limit is constant in $\F$.

\begin{theorem}\label{thm3}
Assume that $W$ satisfies {\rm (i)-(iii)}. 
Then the sequence of functionals $\{\frac{1}{h}\I^{(h)}\}$ 
$\Gamma$-converges, as $h\to 0^{+}$, to the functional

\begin{equation}\label{gammalimit}
\I(y)=
\begin{cases}
\gh     &  \text{ if } F\in\F\,,\\
+\infty   &  \text{ otherwise }\,,
\end{cases}
\end{equation}

\no
with respect to the weak convergence in $L^{p}(\Om;\Mtt)$.
\end{theorem}

\begin{proof}
The proof is devided into two parts. \\

\no
1. {\em Liminf inequality}\\

\no
Let $\{\Fh\}\subset \F_{h}$, 
$\bar{F}\in L^{p}(\Om;\Mtt)$, and let 
$\Fh\weak \bar{F}$ weakly in $L^{p}(\Om;\Mtt)$. 
We have to prove that 

\begin{equation}\label{ineq1}
\I(\bar{F})\leq\liminf_{h\to 0}\frac{1}{h}\I^{(h)}(\Fh) \,.
\end{equation}

\no
We may assume that 
$\liminf_{h\to 0}\frac{1}{h}\I^{(h)}(\Fh)<C$. 
Then, by Theorem \ref{thm1}, there exists 
$F\in L^{p}(\Om;\Mtt)$ 
independent of $x_{2}$ and $x_{3}$, satisfying \eqref{charact}, 
such that $\{\Fhh\}$, up to subsequences, converges to $F$ weakly 
in $L^{p}(\Om;\Mtt)$. 
Hence $\bar{F}^{1}=F^{1}$, $\bar{F}^{2}=\bar{F}^{3}=0$, 
and \eqref{ineq1} hold. 
Moreover, condition \eqref{charact} implies that 
$|\bar{F}^{1}|\leq 1$ a.e. in $(-L,0)$, 
$|\bar{F}^{1}|\leq\zeta_{1}$ a.e. in $(0,L)$, 
and therefore $\bar{F}\in\F$. \\

\no
2. {\em Limsup inequality}\\

\no
We have to show that for each $F\in L^{p}(\Om;\Mtt)$ there exists a sequence $\{\Fh\}\subset\F^{(h)}$ such that $\Fh\weak F$  weakly in 
$L^{p}(\Om;\Mtt)$ and 

\begin{equation}\label{ineq2}
\limsup_{h\to 0}\frac{1}{h}\I^{(h)}(\Fh)\leq \I(F) \,.
\end{equation}

\no
We can assume that $F\in\F$. 
One can construct a pair of measurable functions  
$d_{2},d_{3}\in L^{\infty}((-L,L);\R^{3})$ such that 

\begin{equation*}
(F^{1},d_{2},d_{3})\in \co(SO(3)) \text{ in }(-L,0)\times S\,,
\qquad 
(F^{1},d_{2},d_{3})\in \co(SO(3)H) \text{ in }(0,L)\times S\,.
\end{equation*}

\no
A more detailed construction is contained in \cite{mm}.
We now apply Theorem \ref{thm2} to find a sequence 
$\{\Fh\}\subset\F_{h}$ such that $\Fhh$ converges to 
$(F^{1},d_{2},d_{3})$ weakly in $L^{p}(\Om;\Mtt)$, which 
implies that $\Fh$ converges to $(F^{1},0,0)$. Moreover \eqref{ineq2} holds.

\end{proof}

\medskip
\centerline{\sc Acknowledgements}
\vspace{2mm}

\no
This work was started when MP was post-doc at the Max Planck Institute 
for Mathematics in the Sciences in Leipzig, and carried out when she was 
post-doc in SISSA (Trieste).
The work of MP was partially supported by the DFG Forschergruppe 
``Architektur von nano- und mikrodimensionalen Strukturelementen'' (FOR 522). 


\end{document}